# ITERATIVE ESTIMATING EQUATIONS: LINEAR CONVERGENCE AND ASYMPTOTIC PROPERTIES

By Jiming Jiang,[1] Yihui Luan[2] and You-Gan Wang

*University of California, Davis, Shandong University and CSIRO Mathematical and Information Sciences*

We propose an iterative estimating equations procedure for analysis of longitudinal data. We show that, under very mild conditions, the probability that the procedure converges at an exponential rate tends to one as the sample size increases to infinity. Furthermore, we show that the limiting estimator is consistent and asymptotically efficient, as expected. The method applies to semiparametric regression models with unspecified covariances among the observations. In the special case of linear models, the procedure reduces to iterative reweighted least squares. Finite sample performance of the procedure is studied by simulations, and compared with other methods. A numerical example from a medical study is considered to illustrate the application of the method.

**1. Introduction.** Longitudinal data is often encountered in medical research and economics studies. In the analysis of longitudinal data (e.g., Diggle, Liang and Zeger [3]), the problem of main interest is often related to the estimation of the mean responses which, under a suitable parametric or semiparametric model, depend on a vector $\beta$ of unknown parameters. However, the mission is complicated by the fact that the responses are correlated and the correlations are unknown.

Suppose that $Y$ is a vector of responses that is associated with a matrix $X$ of covariates, which may also be random. Suppose that the (conditional) mean of $Y$ is associated with a vector of parameters, $\theta$. For notational

Received November 2005; revised January 2007.

[1]Supported in part by NSF Grants SES-99-78101 and DMS-04-02824. Part of this work was done while the first author was visiting the Institute for Mathematical Sciences, National University of Singapore in 2005. The visit was supported by the Institute.

[2]Supported in part by NSF of China Grant 10441004.

simplicity, write $\mu = \mu(X,\theta) = \mathrm{E}_\theta(Y|X)$ and $V = \mathrm{Var}(Y|X)$. Hereafter Var or E without the subscript $\theta$ is meant to be taken at the true $\theta$. Consider the following class of estimating functions $\mathcal{G} = \{G = A(Y - \mu)\}$, where $A = A(X,\theta)$. Suppose that $V$ is known. Then, by Theorem 2.1 of Heyde [7], it is easy to show that the optimal estimating function within $\mathcal{G}$ is given by $G^* = \dot\mu' V^{-1}(Y - \mu)$, that is, with $A^* = \dot\mu' V^{-1}$. Therefore, the optimal estimating equation is given by $\dot\mu' V^{-1}(Y - \mu) = 0$. In the case of longitudinal data, the responses are clustered according to the subjects. Let $Y_i$ be the vector of responses collected from the $i$th subject, and $X_i$ the matrix of covariates associated with $Y_i$. Let $\mu_i = \mathrm{E}(Y_i|X_i) = \mu_i(X_i, \beta)$, where $\beta$ is a vector of unknown parameters. Then, under the assumption that $(X_i, Y_i)$, $i = 1, \ldots, n$, are uncorrelated with known $V_i = \mathrm{Var}(Y_i|X_i)$, the optimal estimating equation (for $\beta$) is given by $\sum_{i=1}^n \dot\mu_i' V_i^{-1}(Y_i - \mu_i) = 0$, which is known as the generalized estimating equations (GEE) (e.g., Diggle, Liang and Zeger [3]).

However, the optimal GEE depends on $V_i$, $1 \le i \le n$, which are usually unknown in practice. Therefore, the optimal GEE estimator is not computable. The problem of obtaining an (asymptotically) optimal, or efficient, estimator of $\beta$ without knowing the $V_i$'s is the main concern of the current paper. To motivate our approach, let us first consider a simple example. Suppose that the $Y_i$'s satisfy a (classical) linear model, that is, $\mathrm{E}(Y_i) = X_i\beta$, where $X_i$ is a matrix of fixed covariates, $\mathrm{Var}(Y_i) = V_1$, where $V_1$ is an unknown fixed covariance matrix, and $Y_1, \ldots, Y_n$ are independent. If $V_1$ is known, $\beta$ can be estimated by the following best linear unbiased estimator (BLUE), which is the optimal GEE estimator in this special case:

$$(1.1) \qquad \hat\beta_{\mathrm{BLUE}} = \left(\sum_{i=1}^n X_i' V_1^{-1} X_i\right)^{-1} \sum_{i=1}^n X_i' V_1^{-1} Y_i,$$

provided that $\sum_{i=1}^n X_i' V_1^{-1} X_i$ is nonsingular. On the other hand, if $\beta$ is known, $V_1$ can be estimated consistently as

$$(1.2) \qquad \hat V_1 = \frac{1}{n} \sum_{i=1}^n (Y_i - X_i\beta)(Y_i - X_i\beta)'.$$

It is clear that there is a cycle, which motivates the following algorithm, when neither $\beta$ nor $V_1$ is known. Starting with the identity matrix $I$ for $V_1$, use (1.1) with $V_1$ replaced by $I$ to obtain the initial estimator for $\beta$, which is known as the ordinary least squares (OLS) estimator; then use (1.2) with $\beta$ replaced by the OLS estimator to update $V_1$; then use (1.1) with the new $V_1$ to update $\beta$, and so on. The procedure is expected to result in an estimator that is, at least, more efficient than the OLS estimator; but can we ask a further question, that is, is the procedure going to produce an estimator that is asymptotically as efficient as the BLUE? Before this question can be answered, however, another issue needs to be resolved first, that is, does the iterative procedure converge? These questions will be answered in the sequel.

The example considered above is called balanced data, in which the observations are collected at a common set of times for all the subjects. In this case, the procedure (without iterations) is known as robust estimation for analysis of longitudinal data. However, as pointed out by Diggle, Liang and Zeger [3], page 77, so far the method has been restricted to the case of balanced data. In many practical situations, however, the data is seriously unbalanced, that is, the time sets at which the observations are collected vary from subject to subject. An extension of the robust estimation method to unbalanced longitudinal data will be given.

In fact, we are going to consider a statistical model for the longitudinal data that is much more general than the above example, or the case of balanced data, and propose a similar iterative procedure that not only converges, but converges exponentially. The latter is the main theoretical finding of this paper. Furthermore, at convergence the iterative procedure produces an estimator of $\beta$ that is asymptotically as efficient as the optimal GEE estimator. In Section 2 we propose a semiparametric regression model for the longitudinal data. An iterative estimating equations (IEE) procedure is proposed in Section 3 under the semiparametric model. The convergence property as well as asymptotic behavior of the limiting estimator are studied in Section 4. In Section 5 we consider a special case of our general model, the case of (classical) linear models. Some simulation studies are carried out in Section 6 to investigate empirically the performance of the IEE and its comparison with other methods. A numerical example using data from a medical study is considered in Section 7. Some further discussion and concluding remarks are given in Section 8. Proofs and other technical details are deferred to Section 9.

**2. A semiparametric regression model.** We consider a follow-up study conducted over a set of prespecified visit times $t_1, \ldots, t_b$. Suppose that the responses are collected from subject $i$ at the visit times $t_j$, $j \in J_i \subset J = \{1, \ldots, b\}$. Let $Y_i = (Y_{ij})_{j \in J_i}$. Here we allow the visit times to be dependent on the subject. Let $X_{ij} = (X_{ijl})_{1 \leq l \leq p}$ represent a vector of explanatory variables associated with $Y_{ij}$ so that the first component of $X_{ij}$ corresponds to an intercept (i.e., $X_{ij1} = 1$). Write $X_i = (X_{ij})_{j \in J_i} = (X_{ijl})_{j \in J_i, 1 \leq l \leq p}$. Note that $X_i$ may include both time dependent and independent covariates so that, w.l.o.g., it may be expressed as $X_i = (X_{i1}, X_{i2})$, where the elements of $X_{i1}$ do not depend on $j$ (i.e., time) while those of $X_{i2}$ do. We assume that $(X_i, Y_i)$, $i = 1, \ldots, n$, are independent. Furthermore, it is assumed that

$$\text{E}(Y_{ij}|X_i) = g_j(X_i, \beta), \tag{2.1}$$

where $\beta$ is a $p \times 1$ vector of unknown regression coefficients and $g_j(\cdot, \cdot)$ are fixed functions. We use the notation $\mu_{ij} = \mathrm{E}(Y_{ij}|X_i)$ and $\mu_i = (\mu_{ij})_{j \in J_i}$ throughout this paper. Note that $\mu_i = \mathrm{E}(Y_i|X_i)$. In addition, denote the (conditional) covariance matrix of $Y_i$ given $X_i$ as

$$(2.2) \qquad V_i = \mathrm{Var}(Y_i|X_i),$$

whose $(j,k)$th element is $v_{ijk} = \mathrm{cov}(Y_{ij}, Y_{ik}|X_i) = \mathrm{E}\{(Y_{ij} - \mu_{ij})(Y_{ik} - \mu_{ik})|X_i\}$, $j, k \in J_i$. Note that the dimension of $V_i$ may depend on $i$. Let $D = \{(j,k): j, k \in J_i \text{ for some } 1 \leq i \leq n\}$. The following assumption is made.

ASSUMPTION A1. For any $(j,k) \in D$, the number of different $v_{ijk}$'s is bounded, that is, for each $(j,k) \in D$, there is a set of numbers $\mathcal{V}_{jk} = \{v(j,k,l), 1 \leq l \leq L_{jk}\}$, where $L_{jk}$ is bounded, such that $v_{ijk} \in \mathcal{V}_{jk}$ for any $1 \leq i \leq n$ with $j, k \in J_i$.

Since the number of visit times, $b$, is assumed fixed, an equivalence to Assumption A1 follows.

LEMMA 1. *Assumption* A1 *holds if and only if the number of different $V_i$'s is bounded, that is, there is a set of covariance matrices $\mathcal{V} = \{V(1), \ldots, V(L)\}$, where $L$ is bounded, such that $V_i \in \mathcal{V}$ for any $1 \leq i \leq n$.*

As will be seen, Assumption A1 is essential for consistent estimation of the unknown covariances. On the other hand, the assumption is not unreasonable. We consider some examples.

EXAMPLE 1 (*Balanced data*). Consider the case of balanced data, that is, $J_i = J$, $1 \leq i \leq n$. Furthermore, suppose that $v_{ijk}$ is unknown, but depends only on $j$ and $k$. Then, one has $V_i = V_1$, an unknown constant covariance matrix, $1 \leq i \leq n$. Thus, Assumption A1 is satisfied.

The following examples show that Assumption A1 remains valid in some unbalanced cases as well.

EXAMPLE 2. Suppose that the observational times are equally spaced. In such a case we may assume, without loss of generality, that $t_j = j$. Suppose that the responses $Y_{ij}$ satisfy

$$(2.3) \qquad Y_{ij} = x'_{ij}\beta + u_i + w_{ij} + e_{ij},$$

$i = 1, \ldots, n$, $j \in J_i \subset J$, where $x_{ij}$ is a vector of *fixed* covariates, $u_i$ is a subject-specific random effect, $w_{ij}$ corresponds to a serial correlation, and $e_{ij}$ represents a measurement error. It is assumed that the $u_i$'s are independent and distributed as $N(0, \sigma_u^2)$, and the $e_{ij}$'s are independent and distributed as $N(0, \sigma_e^2)$. As for the $w_{ij}$'s, it is assumed that they satisfy the equation of the first-order autoregressive, or AR(1), process, $w_{ij} = \phi w_{ij-1} + z_{ij}$, where $\phi$ is a constant and $|\phi| < 1$, and the $z_{ij}$'s are independent with distribution $N\{0, \sigma_w^2(1-\phi^2)\}$. Finally, we assume that $u$, $w$ and $e$ are independent. Under this model, we have $\mathrm{E}(Y_i|X_i) = X_i\beta$, where $X_i = (x'_{ij})_{j \in J_i}$. Furthermore, we have $V_i = (v_{ijk})_{j,k \in J_i}$, where $v_{ijk} = \sigma_u^2 + \sigma_w^2 \phi^{|j-k|} + \sigma_e^2 \delta_{j,k}$, and $\delta_{j,k} = 1$ if $j = k$ and 0 otherwise (e.g., Anderson [1], page 174). It follows that $v_{ijk}$ does not depend on $i$. We now further specify the $J_i$'s. Suppose that the responses are collected over a week such that the data is collected on Monday, Wednesday and Friday for 40 percent of the subjects. For the rest of the subjects, the data is collected on Tuesday and Thursday. Then the $J_i$'s are either $\{1,3,5\}$ or $\{2,4\}$. It follows that $D$ consists of $(j,j)$, $j = 1,\ldots,5$, $(1,3)$, $(1,5)$, $(3,5)$, $(2,4)$ and the pairs with the components exchanged. Furthermore, the $V_i$'s are either $V^{(1)}$ or $V^{(2)}$, where

$$V^{(1)} = \sigma_u^2 \begin{pmatrix} 1 & 1 \\ 1 & 1 \end{pmatrix} + \sigma_w^2 \begin{pmatrix} 1 & \phi^2 \\ \phi^2 & 1 \end{pmatrix} + \sigma_e^2 \begin{pmatrix} 1 & 0 \\ 0 & 1 \end{pmatrix},$$

$$V^{(2)} = \sigma_u^2 \begin{pmatrix} 1 & 1 & 1 \\ 1 & 1 & 1 \\ 1 & 1 & 1 \end{pmatrix} + \sigma_w^2 \begin{pmatrix} 1 & \phi^2 & \phi^4 \\ \phi^2 & 1 & \phi^2 \\ \phi^4 & \phi^2 & 1 \end{pmatrix} + \sigma_e^2 \begin{pmatrix} 1 & 0 & 0 \\ 0 & 1 & 0 \\ 0 & 0 & 1 \end{pmatrix}.$$

In practice, the covariance structure of the data may not be known. For example, if the $w_{ij}$'s are known to satisfy a process with covariances $\mathrm{cov}(w_{ij}, w_{ik}) = \sigma_w^2 \rho(j,k)$ with $\rho(\cdot,\cdot)$ completely unknown, then the covariance matrices $V^{(1)}$ and $V^{(2)}$ will be practically unspecified. However, there are only two possible covariance matrices. Thus, again by Lemma 1, Assumption A1 is satisfied.

EXAMPLE 3 (*Growth curve*). This model is the same as Example 2 except that a term $(\mu + a_i)t_{ij}$ is added to the right-hand side of (2.3), where $\mu$ is an unknown baseline slope, $a_i$ is a random effect with distribution $N(0, \sigma_a^2)$, and $t_{ij}$ is the time at which response $Y_{ij}$ is collected. Assume that $a_1, \ldots, a_n$ are independent, and the $a$'s are independent with $u$, $w$ and $e$. Then we have $v_{ijk} = \sigma_u^2 + \sigma_a^2 t_{ij} t_{ik} + \sigma_w^2 \phi^{|j-k|} + \sigma_e^2 \delta_{j,k}$. Note that, unlike Example 2, here the expression of $v_{ijk}$, indeed, depends on $i$, the subject. However, since $t_{ij} \in J = \{1, \ldots, b\}$, the number of different values of the product $t_{ij}t_{ik}$ is bounded. As a result, the number of different $v_{ijk}$'s is bounded, and hence Assumption A1 is satisfied.

As in the previous example, the parametric covariance structure may not be known in practice, or one may be concerned with model misspecifications. Therefore, as a robust approach, one may prefer to use unspecified covariances with the understanding that Assumption A1 is satisfied.

In a way, the models in the previous examples are formulated as the classical linear models, in which the covariates $x_{ij}$ are considered fixed. We now consider a logistic model in which the covariates are considered random.

EXAMPLE 4. Suppose that $Y_{ij}$, $i = 1, \ldots, n$, $j \in J_i$, are binary responses (i.e., $Y_{ij} = 0$ or 1). A covariate is associated with the age of the subject. More specifically, there are nine age groups: 40–44, ..., 75–79 and 80 or over. Let $X_i = j$ if the age of subject $i$ belongs to the $j$th group. Suppose that, given $X_i$ and an additional (unobservable) subject-specific random effect $\alpha_i$, $Y_{ij}$, $j \in J_i$, are (conditionally) independent such that $\text{logit}\{P(Y_{ij} = 1|X_i, \alpha_i)\} = \beta_0 + \beta_1 X_i + \alpha_i$, where $\text{logit}(x) = \log\{x/(1-x)\}$, and $\beta_0$ and $\beta_1$ are unknown regression coefficients. Furthermore, $(Y_i, X_i, \alpha_i)$, $i = 1, \ldots, n$, are independent, where $Y_i = (Y_{ij})_{j \in J_i}$. Finally, we assume that $\alpha_i|X_i \sim N(0, \sigma^2)$. It follows that $E(Y_{ij}|X_i) = E\{E(Y_{ij}|X_i, \alpha_i)|X_i\} = E\{h(\beta_0 + \beta_1 X_i + \alpha_i)|X_i\} = E\{h(\beta_0 + \beta_1 X_i + \sigma\xi)\}$, where $h(x) = e^x/(1 + e^x)$, and the last expectation is with respect to $\xi \sim N(0, 1)$. Similarly, we have $v_{ijk} = E[\{h(\beta_0 + \beta_1 X_i + \sigma\xi)\}^{2-\delta_{j,k}}] - [E\{h(\beta_0 + \beta_1 X_i + \sigma\xi)\}]^2$. If the conditional independence is not assumed, neither is any other specified conditional covariance structure, and the above expression for $v_{ijk}$ may not hold. In such a case, the $v_{ijk}$'s are unspecified. However, it is easy to see that for any $j, k$, there are, at most, nine different $v_{ijk}$'s (corresponding to $X_i = 1, \ldots, 9$). Thus, Assumption A1 is satisfied.

## 3. Iterative estimating equations.

3.1. *Estimation of $\beta$ when $V_i$'s are known.* Our main interest is to estimate $\beta$, the vector of regression coefficients. According to the earlier discussion, if the $V_i$'s are known, $\beta$ may be estimated by the GEE given below,

$$(3.1) \qquad \sum_{i=1}^{n} \dot{\mu}_i' V_i^{-1}(Y_i - \mu_i) = 0,$$

where $\mu_i = (\mu_{ij})_{j \in J_i}$ with $\mu_{ij} = g_j(X_i, \beta)$, and $\dot{\mu}_i$ is the matrix of first derivatives whose $(j, l)$ element is given by $\partial \mu_{ij}/\partial \beta_l$, $j \in J_i$, $1 \leq l \leq p$ (e.g., Liang and Zeger [12]). Namely, the estimator, $\hat{\beta}$, is defined as the solution to (3.1).

3.2. *Estimation of $V_i$'s when $\beta$ is known.* On the other hand, if $\beta$ is known, the covariance matrices $V_i$ can be estimated by the method of moments (MoM) as follows. Let $I(j, k, l)$ denote the set of indexes $1 \leq i \leq n$ such that $v_{ijk} = v(j, k, l)$ (see Assumption A1 in Section 2). For any $(j, k) \in D$, $1 \leq l \leq L_{jk}$, define

$$(3.2) \quad \hat{v}(j, k, l) = \frac{1}{n(j, k, l)} \sum_{i \in I(j, k, l)} \{Y_{ij} - g_j(X_i, \beta)\}\{Y_{ik} - g_k(X_i, \beta)\},$$

where $n(j,k,l) = |I(j,k,l)|$, the cardinality. Then, define $\hat{V}_i = (\hat{v}_{ijk})_{j,k \in J_i}$, where $\hat{v}_{ijk} = \hat{v}(j,k,l)$, if $i \in I(j,k,l)$. In this approach, the estimators of the covariances are obtained componentwise. Alternatively, the estimators may be obtained as matrices. Let $I(l)$ denote the set of indexes $1 \leq i \leq n$ such that $V_i = V(l)$ (see Lemma 1). For any $1 \leq l \leq L$, define

$$(3.3) \qquad \tilde{V}(l) = \frac{1}{n(l)} \sum_{i \in I(l)} (Y_i - \mu_i)(Y_i - \mu_i)',$$

where $n(l) = |I(l)|$. Then define $\tilde{V}_i = \tilde{V}(l)$ if $i \in I(l)$. The following lemma justifies the use of the words "method of moments." Write $X = (X_i)_{1 \leq i \leq n}$.

LEMMA 2. *Under the assumed model, we have* $\mathrm{E}(\hat{V}_i|X) = V_i = \mathrm{E}(\tilde{V}_i|X)$, $1 \leq i \leq n$.

The proof is straightforward and therefore omitted. In some cases, the two estimators are identical, as the following lemma shows.

LEMMA 3. *Suppose that there is a division* $J = \{1, \ldots, b\} = J(1) \cup \cdots \cup J(s)$ *such that* $J(q) \cap J(r) = \varnothing$, *if* $q \neq r$ *and for each* $1 \leq i \leq n$ *there is* $1 \leq r \leq s$ *such that* $J_i = J(r)$. *Then, we have* $\hat{V}_i = \tilde{V}_i$, $1 \leq i \leq n$.

The proof follows directly from the definitions. Again, we consider some examples.

EXAMPLE 1 (*Continued*). Here $L_{jk} = L = 1$ and $I(j,k,1) = I(1) = \{1, \ldots, n\}$. So that $\hat{v}_{ijk} = \hat{v}(j,k,1) = n^{-1} \sum_{i=1}^{n} \{Y_{ij} - g_j(X_i, \beta)\}\{Y_{ik} - g_k(X_i, \beta)\}$, the $(j,k)$ element of $\tilde{V}_i = \tilde{V}_1 = n^{-1} \sum_{i=1}^{n} (Y_i - \mu_i)(Y_i - \mu_i)'$.

EXAMPLE 2 (*Continued*). Consider the special case in which $J_i$ is either $\{1,3,5\}$ or $\{2,4\}$. Let $I(1) = \{1 \leq i \leq n : J_i = \{1,3,5\}\}$, $I(2) = \{1 \leq i \leq n : J_i = \{2,4\}\}$ and $n(l) = |I(l)|$, $l = 1, 2$. Then we have

$$\hat{v}(j,k,1) = \frac{1}{n(1)} \sum_{J_i = \{1,3,5\}} (Y_{ij} - x'_{ij}\beta)(Y_{ik} - x'_{ik}\beta), \qquad j,k = 1,3,5;$$

$$\hat{v}(j,k,2) = \frac{1}{n(2)} \sum_{J_i = \{2,4\}} (Y_{ij} - x'_{ij}\beta)(Y_{ik} - x'_{ik}\beta), \qquad j,k = 2,4.$$

On the other hand, we have

$$\tilde{V}(1) = \frac{1}{n(1)} \sum_{J_i = \{1,3,5\}} (Y_i - X_i\beta)(Y_i - X_i\beta)',$$

$$\tilde{V}(2) = \frac{1}{n(2)} \sum_{J_i = \{2,4\}} (Y_i - X_i\beta)(Y_i - X_i\beta)',$$

where $X_i = (x'_{ij})_{j \in J_i}$. It follows that $\hat{V}_i = \tilde{V}_i$, $1 \leq i \leq n$.

The next example shows that, when the two estimators are different, the first estimator is likely to be more efficient than the second one.

EXAMPLE 5. Consider, again, Example 2, but this time assume that every subject has a baseline measure on Monday. Then 50% of the subjects have follow-up measures on Tuesday, and 50% of the subjects have the follow-ups on Wednesday. In other words, $J_i$ is either $\{1,2\}$ or $\{1,3\}$. To be more specific, let $n = 2m$ such that subjects $1, \ldots, m$ have the follow-ups on Tuesday, and subjects $m+1, \ldots, n$ have the follow-ups on Wednesday. If we assume that the covariances do not depend on $i$, and write $v_{jk} = \text{cov}(Y_{ij}, Y_{ik})$, then there are two different covariance matrices,

$$V(1) = \begin{pmatrix} v_{11} & v_{12} \\ v_{21} & v_{22} \end{pmatrix}, \qquad V(2) = \begin{pmatrix} v_{11} & v_{13} \\ v_{31} & v_{33} \end{pmatrix}.$$

On the other hand, it is easy to see that $L_{jk} = 1$ for all possible $j$, $k$, and $I(1,1,1) = \{1, \ldots, n\}$, $I(1,2,1) = I(2,2,1) = \{1, \ldots, m\}$ and $I(1,3,1) = I(3,3,1) = \{m+1, \ldots, n\}$. Thus, we have

$$\hat{v}_{11} = \frac{1}{n} \sum_{i=1}^{n} (Y_{i1} - x'_{i1}\beta)^2, \tag{3.4}$$

$\hat{v}_{12} = m^{-1} \sum_{i=1}^{m} (Y_{i1} - x'_{i1}\beta)(Y_{i2} - x'_{i2}\beta)$, $\hat{v}_{22} = m^{-1} \sum_{i=1}^{m} (Y_{i2} - x'_{i2}\beta)^2$, $\hat{v}_{13} = m^{-1} \sum_{i=m+1}^{n} (Y_{i1} - x'_{i1}\beta)(Y_{i3} - x'_{i3}\beta)$ and $\hat{v}_{33} = m^{-1} \sum_{i=m+1}^{n} (Y_{i3} - x'_{i3}\beta)^2$. Similarly, we have

$$\tilde{V}(1) = \frac{1}{m} \sum_{i=1}^{m} (Y_i - X_i\beta)(Y_i - X_i\beta)',$$

$$\tilde{V}(2) = \frac{1}{m} \sum_{i=m+1}^{n} (Y_i - X_i\beta)(Y_i - X_i\beta)'.$$

Comparing the estimators, it is seen that $\hat{v}$ and $\tilde{v}$ are identical except for the case $j = k = 1$. In the latter case, the componentwise estimator of $v_{11}$, $\hat{v}_{11}$, is given by (3.4); on the other hand, there are two estimators of $v_{11}$ by the matrix approach, namely,

$$\tilde{v}_{11}(1) = \frac{1}{m} \sum_{i=1}^{m} (Y_{i1} - x'_{i1}\beta)^2,$$

$$\tilde{v}_{11}(2) = \frac{1}{m} \sum_{i=m+1}^{n} (Y_{i1} - x'_{i1}\beta)^2.$$

However, none of the latter estimators is as efficient as $\hat{v}_{11}$, because each of them is based on half of the samples.

For such a reason, in the sequel we mainly focus on the elementwise estimators.

3.3. *Iterative procedure.* The main points of the previous subsections may be summarized as follows: If the $V_i$'s were known, one could estimate $\beta$ by the GEE; on the other hand, if $\beta$ were known, one could estimate the $V_i$'s by the MoM. It is clear that there is a cycle here, which motivates the following iterative procedure. Starting with an initial estimator of $\beta$, use (3.2), with $\beta$ replaced by the initial estimator, to obtain the estimators of the covariances; then use (3.1) to update the estimator of $\beta$; and repeat the process. We call such a procedure iterative estimating equations, or IEE. If the procedure converges, the limiting estimator is called the IEE estimator.

In practice, the initial estimate of $\beta$ may be obtained as the solution to (3.1) with $V_i = I$, the identity matrix (with the suitable dimension). It can be shown (see Jiang, Luan and Wang [10]) that in the case of balanced data of Example 1, procedure is equivalent to maximum likelihood (ML), if the errors $\varepsilon_i = Y_i - \mu_i$ are normal. Thus, without the normality assumption, IEE may be viewed as a quasi-likelihood method. Although for unbalanced data IEE is not equivalent to ML even under normality, the procedure is expected to improve the efficiency of the estimator of $\beta$. To see this, note that if the initial estimator of $\beta$ is consistent, the estimators of the $V_i$'s by (3.2) are expected to be consistent; and, with consistent estimators of the $V_i$'s, the efficiency of the next step estimator of $\beta$ should improve (because $\hat{V}_i \approx V_i$), and so on. Furthermore, IEE is very easy to operate and is, in fact, no stranger to practitioners. For example, the procedure has some similarity with the so-called backfitting algorithm (e.g., Hastie and Tibshirani [6], Section 4.4). But before we study the efficiency and other potentially nice properties of this procedure, we need to make sure that it converges. The next question is: Given the convergence of IEE, is the limiting estimator asymptotically as efficient as the GEE estimator that one would have obtained if the $V_i$'s were known? These fundamental questions will be answered in the next section.

## 4. Convergence and asymptotic properties.

4.1. *Linear convergence.* First we formulate the IEE procedure as follows. Let $v = (v_r)_{1 \leq r \leq R}$ denote the vector of different covariances involved in the $V_i$'s. We assume that for any $v$, there is a unique solution to (3.1) for $\beta$. Denote this solution by $\beta(v)$. Also, for any $\beta$, let $v(\beta)$ denote the vector $v$ whose corresponding components are given by the right-hand side of (3.2). We use an example to illustrate.

EXAMPLE 5 (*Continued*). Consider, once again, Example 5. In this case, we have $v = (v_{11}, v_{12}, v_{22}, v_{13}, v_{33})'$. Equation (3.1) has a unique solution

given by
$$\beta(v) = \left(\sum_{i=1}^{n} X_i' V_i^{-1} X_i\right)^{-1} \sum_{i=1}^{n} X_i' V_i^{-1} Y_i,$$
where $X_i = (x_{ij}')_{j \in J_i}$, $J_i = \{1, 2\}$, $V_i = V(1)$, $1 \leq i \leq m$, and $J_i = \{1, 3\}$, $V_i = V(2)$, $m+1 \leq i \leq n$. Furthermore, $v(\beta)$ is the vector whose components are given by the right-hand sides of (3.4) and the four equations below.

With such notation, the IEE may be formulated as follows: Take the initial $v$ as $\hat{v}^{(0)}$ whose component is 0, if it is a covariance (i.e., $j \neq k$); otherwise, the component is 1. This means that, initially, all the $V_i$'s are taken as the identity matrices of suitable dimensions. Then we have $\hat{\beta}^{(1)} = \beta\{\hat{v}^{(0)}\}$, $\hat{v}^{(1)} = \hat{v}^{(0)}$; $\hat{\beta}^{(2)} = \hat{\beta}^{(1)}$, $\hat{v}^{(2)} = v\{\hat{\beta}^{(1)}\}$, and so on. In general, we have $\hat{\beta}^{(2m-1)} = \beta\{\hat{v}^{(2m-2)}\}$, $\hat{v}^{(2m-1)} = \hat{v}^{(2m-2)}$; $\hat{\beta}^{(2m)} = \hat{\beta}^{(2m-1)}$, $\hat{v}^{(2m)} = v\{\hat{\beta}^{(2m-1)}\}$ for $m = 1, 2, \ldots$. There is an alternative expression, which is more convenient in terms of obtaining the convergence rate. Write $\tau = \beta \circ v$ and $\rho = v \circ \beta$. Then we have $\hat{\beta}^{(2m)} = \hat{\beta}^{(2m-1)}$, $\hat{\beta}^{(2m+1)} = \tau\{\hat{\beta}^{(2m-1)}\}$, $\hat{v}^{(2m-1)} = \hat{v}^{(2m-2)}$, $\hat{v}^{(2m)} = \rho\{\hat{v}^{(2m-2)}\}$ for $m = 1, 2, \ldots$. Theorem 1 below states that, under very mild conditions, the IEE converges linearly with probability tending to one. Here we adapt a term from numerical analysis: An iterative algorithm that results in a sequence $x^{(m)}$, $m = 1, 2, \ldots$, converges linearly to a limit $x^*$, if there is $0 < \varrho < 1$ such that $\sup_{m \geq 1}\{|x^{(m)} - x^*|/\varrho^m\} < \infty$ (e.g., Press et al. [14]).

Let $L_1 = \max_{1 \leq i \leq n} \max_{j \in J_i} s_{ij}$ with $s_{ij} = \sup_{|\tilde{\beta} - \beta| \leq \varepsilon_1} |(\partial/\partial \beta) g_j(X_i, \tilde{\beta})|$, where $\beta$ represents the true parameter vector, $\varepsilon_1$ is any positive constant, and $(\partial/\partial \beta) f(\tilde{\beta})$ means $(\partial f/\partial \beta)|_{\beta = \tilde{\beta}}$. Similarly, let $L_2 = \max_{1 \leq i \leq n} \max_{j \in J_i} w_{ij}$, where $w_{ij} = \sup_{|\tilde{\beta} - \beta| \leq \varepsilon_1} \|(\partial^2/\partial \beta \, \partial \beta') g_j(X_i, \tilde{\beta})\|$ ($\|\cdot\|$ is the spectral norm defined in Section 9). Furthermore, define the set $\mathcal{V} = \{v : \lambda_{\min}(V_i) \geq \delta_0, \lambda_{\max}(V_i) \leq M_0, 1 \leq i \leq n\}$, where $\lambda_{\min}$ and $\lambda_{\max}$ represent the smallest and largest eigenvalues, respectively, and $\delta_0$ and $M_0$ are given positive constants. Note that $\mathcal{V}$ is a nonrandom set.

An array of nonnegative definite matrices $\{A_{n,i}\}$ is bounded from above if $\|A_{n,i}\| \leq c$ for some constant $c$; the array is bounded from below if $A_{n,i}^{-1}$ exists and $\|A_{n,i}^{-1}\| \leq c$ for some constant $c$. A sequence of random matrices is bounded in probability, denoted by $A_n = O_P(1)$, if for any $\varepsilon > 0$, there are $M > 0$ and $N \geq 1$ such that $P(\|A_n\| \leq M) > 1 - \varepsilon$, if $n \geq N$. The sequence is bounded away from zero in probability if $A_n^{-1} = O_P(1)$. Note that the definition also applies to a sequence of random variables, considered as a special case of random matrices.

Recall that $p$ is the dimension of $\beta$ and $R$ the dimension of $v$. Also recall assumption Assumption A1 (above Lemma 1). We make the following additional assumptions.

ASSUMPTION A2. The functions $g_j(X_i, \beta)$ are twice continuously differentiable with respect to $\beta$; $\mathrm{E}(|Y_i|^4)$, $1 \leq i \leq n$, are bounded; and $L_1$, $L_2$, $\max_{1 \leq i \leq n}(\|V_i\| \vee \|V_i^{-1}\|)$ are $O_\mathrm{P}(1)$.

ASSUMPTION A3 (*Consistency of GEE estimator*). For any given $V_i$, $1 \leq i \leq n$, bounded from above and below, the GEE equation (3.1) has a unique solution $\hat{\beta}$ that is consistent.

ASSUMPTION A4 (*Differentiability of GEE solution*). For any $v$, the solution to (3.1), $\beta(v)$, is continuously differentiable with respect to $v$, and $\sup_{v \in \mathcal{V}} \|\partial \beta / \partial v\| = O_\mathrm{P}(1)$.

ASSUMPTION A5. $n(j, k, l) \to \infty$ for any $1 \leq l \leq L_{jk}$, $(j, k) \in D$, as $n \to \infty$.

THEOREM 1. *Under Assumptions* A1–A5, $\mathrm{P}(IEE\ converges) \to 1$ *as* $n \to \infty$. *Furthermore, we have* $\mathrm{P}[\sup_{m \geq 1}\{|\hat{\beta}^{(m)} - \hat{\beta}^*|/(p\eta)^{m/2}\} < \infty] \to 1$, $\mathrm{P}[\sup_{m \geq 1}\{|\hat{v}^{(m)} - \hat{v}^*|/(R\eta)^{m/2}\} < \infty] \to 1$ *as* $n \to \infty$ *for any* $0 < \eta < (p \vee R)^{-1}$, *where* $(\hat{\beta}^*, \hat{v}^*)$ *is the (limiting) IEE estimator.*

REMARK 1. Because our approach is based on estimating functions (EFs), it is therefore necessary to impose some "regularity" conditions to make sure that the EFs are valid. Assumption A3 may be regarded as such regularity conditions. Note that here we did not specify the technical conditions for the existence, uniqueness and consistency of the solution to the GEE equation (3.1). However, standard conditions for the existence and consistency of the GEE solution have been discussed extensively in the literature. See, for example, Liang and Zeger [12] and Heyde [7], Section 12.2. Furthermore, an easy-to-verify sufficient condition can be given for the uniqueness of the GEE solution, with the additional assumption that the initial estimator, $\hat{\beta}^{(1)}$, is consistent. To see this, note that by the proof of Theorem 1 given in Section 9, it is seen that the linear convergent property of the sequences $\hat{\beta}^{(m)}$ and $\hat{v}^{(m)}$ ($m = 1, 2, \ldots$) only depends on a small neighborhood of the true $\beta$ where the first sequence is confined, if the initial estimator $\hat{\beta}^{(1)}$ is consistent. On the other hand, note that, for fixed $V_i$'s, the left-hand side of (3.1) is the gradient of a function of $\beta$. Namely, write $q(\beta) = \sum_{i=1}^n (Y_i - \mu_i)' V_i^{-1}(Y_i - \mu_i)$. Then, $\partial q / \partial \beta$ is $-2$ times the left-hand side of (3.1). Therefore, any conditions that guarantee (local) convexity of $q(\cdot)$ will be sufficient for the (local) uniqueness of the solution to (3.1) [note that, because of the earlier remark, all one needs is that (3.1) has a unique solution *locally*, i.e., in a small neighborhood of the true $\beta$]. It is easy to show that $\partial^2 q / \partial \beta \, \partial \beta' = 2(I_1 - I_2)$ with $I_1 = \sum_{i=1}^n (\partial \mu_i' / \partial \beta) V_i^{-1} (\partial \mu_i / \partial \beta')$ and $I_2 = \sum_{i=1}^n I_{2,i}$, where the $(k, l)$ element of $I_{2,i}$ is $(\partial^2 \mu_i' / \partial \beta_k \, \partial \beta_l) V_i^{-1}(Y_i - \mu_i)$. Since $I_1$ is nonnegative definite, a condition of (asymptotic) nonsingularity is:

(i) $\liminf_{n\to\infty} \lambda_{\min}(n^{-1}I_1) > 0$.

Furthermore, it is seen that at the true $\beta$, $I_2$ has mean zero. Also, Assumption A2 and boundedness from below of the $V_i$'s imply that $\text{var}(I_2) = O(n)$. Hence, we have $n^{-1}I_2 \to 0$ in $L_2$. Thus, if $I_2$ is uniformly continuous in $\beta$, the matrix $\partial^2 q/\partial\beta\,\partial\beta'$ is expected to be (asymptotically) positive definite in a small neighborhood of the true $\beta$, which implies asymptotic local convexity of $q(\cdot)$. A condition for the (asymptotic) uniform continuity is:

(i) for any $\delta > 0$, there is $\varepsilon > 0$ such that $|\tilde{\mu}_i - \mu_i| < \delta$, $|\tilde{s}_{ikl} - s_{ikl}| < \delta$ for all $1 \leq i \leq n$ and $k, l \in J_i$, if $|\tilde{\beta} - \beta| < \varepsilon$.

Here $s_{ikl}$ denotes $\partial^2 \mu_i/\partial\beta_k\,\partial\beta_l$ evaluated at $\beta$ and $\tilde{s}_{ikl}$ the same quantity evaluated at $\tilde{\beta}$. It can be shown that conditions (i) and (ii) are, indeed, sufficient for the asymptotic uniqueness of the solution to the GEE equation (3.1) (see Jiang, Luan and Wang [10]).

REMARK 2. For the most part, Assumption A4 is about differentiability of the implicit function $\beta(v)$. It simply requires that the $\beta$ estimates not be too sensitive to the changes in $v$. This is quite reasonable because, otherwise, the estimates of $\beta$ will be very unstable so that convergence cannot be achieved. Sufficient conditions for such differentiability can be found in standard texts of advanced calculus. For example, the following result can be implied directly from the theorem on page 265 of Courant and John [2]: Let $\phi$ denote the left-hand side of (3.1) as a (vector-valued) function of $\beta$ and $v$. If equation (3.1) is satisfied at a point $(\tilde{\beta}', \tilde{v}')'$, and the Jacobian of $\phi$ with respect to $\beta$ differs from zero at that point, then in the neighborhood of that point equation (3.1) can be solved in one and only one way for $\beta$, and this solution gives $\beta$ as a continuously differentiable function of $v$. Note that the existence of $(\tilde{\beta}', \tilde{v}')'$ is a consequence of Assumption A3 (see the remark above). Also, as noted above, for the results of Theorem 1 to hold, one only needs the assumptions to be valid in a small neighborhood of $\beta$. Finally, the nonzero-Jacobian property is only required with probability tending to one.

REMARK 3. It is clear that the restriction $\eta < (p \vee R)^{-1}$ is unnecessary [because, for example, $(p\eta_1)^{-m/2} < (p\eta_2)^{-m/2}$ for any $\eta_1 \geq (p \vee R)^{-1} > \eta_2$], but linear convergence would only make sense when $\varrho < 1$ (see the definition above).

REMARK 4. The proof of Theorem 1 in fact demonstrates that for any $\delta > 0$, there are positive constants $M_1$, $M_2$ and integer $N$ that depend only on $\delta$ such that, for all $n \geq N$, we have $\text{P}[\sup_{m\geq 1}\{|\hat{\beta}^{(m)} - \hat{\beta}^*|/(p\eta)^{m/2}\} \leq M_1] > 1 - \delta$, $\text{P}[\sup_{m\geq 1}\{|\hat{v}^{(m)} - \hat{v}^*|/(R\eta)^{m/2}\} \leq M_2] > 1 - \delta$.

4.2. *Asymptotic behavior of the IEE estimator.* In the discussion in Section 3.3 we conjectured that the (limiting) IEE estimator is an efficient estimator. In this subsection we show that this conjecture is indeed true. Since efficiency is usually defined from an asymptotic point of view (e.g., Lehmann [11]), we need to study asymptotic properties of the IEE estimator. The first result is regarding its consistency.

THEOREM 2. *Under the assumptions of Theorem 1, $(\hat{\beta}^*, \hat{v}^*)$ is consistent.*

To establish the asymptotic efficiency, we need to strengthen Assumptions A2 and A5. Define the following: $L_{2,0} = \max_{1 \leq i \leq n} \max_{j \in J_i} \|\partial^2 \mu_{ij}/\partial\beta\,\partial\beta'\|$, $L_3 = \max_{1 \leq i \leq n} \max_{j \in J_i} d_{ij}$, where

$$d_{ij} = \max_{1 \leq a,b,c \leq p} \sup_{|\tilde{\beta}-\beta| \leq \varepsilon_1} \left| \frac{\partial^3}{\partial\beta_a\,\partial\beta_b\,\partial\beta_c} g_j(X_i, \tilde{\beta}) \right|.$$

ASSUMPTION A2'. Same as Assumption A2 except that $g_j(X_i, \beta)$ are three-times continuously differentiable with respect to $\beta$, and that $L_2 = O_P(1)$ is replaced by $L_{2,0} \vee L_3 = O_P(1)$.

ASSUMPTION A5'. There is a positive integer $\gamma$ such that $n/\{n(j,k,l)\}^\gamma \to 0$ for any $1 \leq l \leq L_{jk}$, $(j,k) \in D$, as $n \to \infty$.

We also need the following additional assumption.

ASSUMPTION A6. $n^{-1} \sum_{i=1}^n \dot{\mu}_i' V_i^{-1} \dot{\mu}_i$ is bounded away from zero in probability.

Also, let $\tilde{\beta}$ be the solution to (3.1) with the true $V_i$'s. Note that $\tilde{\beta}$ is efficient, or optimal in the sense discussed in Section 1, but not computable unless the true $V_i$'s are known.

THEOREM 3. *Under Assumptions A1, A2', A3, A4, A5' and A6, we have $\sqrt{n}(\hat{\beta}^* - \tilde{\beta}) \to 0$ in probability. Therefore, asymptotically, $\hat{\beta}^*$ is as efficient as $\tilde{\beta}$.*

The proofs of Theorem 2 and Theorem 3 follow, for the most part, the standard arguments for asymptotic properties of Z-estimators (e.g., van der Vaart and Wellner [16], Section 3.3) and optimality of estimating equations (Godambe [4]), and therefore are omitted. Interested readers are referred to a technical report (Jiang, Luan and Wang [10]). The proof of Theorem 3 also reveals the asymptotic expansion

$$\hat{\beta}^* - \beta = \left(\sum_{i=1}^n \dot{\mu}_i' V_i^{-1} \dot{\mu}_i\right)^{-1} \sum_{i=1}^n \dot{\mu}_i' V_i^{-1}(Y_i - \mu_i) + \frac{o_{\mathrm{P}}(1)}{\sqrt{n}}, \quad (4.1)$$

where $o_{\mathrm{P}}(1)$ represents a term that converges to zero (vector) in probability. By Theorem 3, (4.1) also holds with $\hat{\beta}^*$ replaced by $\tilde{\beta}$, even through the latter is typically not computable. In the next section, we will be looking at a case with an "exact" version of (4.1), that is, the equation without the term $o_{\mathrm{P}}(1)/\sqrt{n}$, for $\tilde{\beta}$.

Suppose that the initial estimator of $\beta$ is consistent. Then, under some regularity conditions, the next-step estimator of $v$ is consistent; furthermore, the next-step estimator of $\beta$, say, $\hat{\beta}^{(3)}$ [according to the notation below Example 5 (continued) in Section 4.1], is not only consistent but also asymptotically as efficient as $\tilde{\beta}$. In other words, one has obtained an efficient estimator of $\beta$ after one iteration of the IEE. Although both $\hat{\beta}^*$ and $\hat{\beta}^{(3)}$ are efficient estimators, there are reasons that $\hat{\beta}^*$ is more desirable in some situations. For example, in the case of balanced data of Example 1 with normality, $\hat{\beta}^*$ is the same as the maximum likelihood estimator (MLE), while $\hat{\beta}^{(3)}$ is not. In Section 6.2 we show by a simple simulated example that IEE may result in substantial improvement over $\hat{\beta}^{(3)}$ in a small sample situation. Although, in general, $\hat{\beta}^*$ may not be the MLE (especially if normality fails), there is a tendency of seeking for continuing improvement, and therefore not to stop after one iteration. It should be pointed out that the continuing improvement may not be as significant as in the early stage when the procedure is close to convergence. Nevertheless, because of the fast convergence (Theorem 1; also see Sections 6 and 7), one may not need to run the IEE for a lot of steps anyway in order to achieve convergence. This, combined with the fact that one has an asymptotically efficient estimator even after one iteration, would give a practitioner double confidence that he/she does not have to run the IEE for many steps in order to use it.

**5. The special case of linear models.** One special case of the semiparametric regression model is the following linear model for longitudinal data. Let $X_i$ be a matrix of fixed covariates associated with the $i$th subject. Suppose that $Y_1,\ldots,Y_n$ are independent such that $\mathrm{E}(Y_i) = X_i\beta$ and $\mathrm{Var}(Y_i) = V_i$, where $V_i$ is an unknown fixed covariance matrix, $1 \leq i \leq n$. Under this model, the GEE (3.1) reduces to $\sum_{i=1}^n X_i'V_i^{-1}(Y_i - X_i\beta) = 0$, which has an explicit solution

$$\hat{\beta} = \left(\sum_{i=1}^n X_i'V_i^{-1}X_i\right)^{-1} \sum_{i=1}^n X_i'V_i^{-1}Y_i, \quad (5.1)$$

provided that the matrix $\sum_{i=1}^n X_i' V_i^{-1} X_i$ is nonsingular. Note that (5.1) is the optimal weighted least squares (WLS) estimator, provided that the $V_i$'s are known. Here a WLS estimator is defined as the solution to the minimization problem $\min_\beta \{(Y - X\beta)' W (Y - X\beta)\}$, where $Y = (Y_i)_{1 \leq i \leq n}$, $X = (X_i)_{1 \leq i \leq n}$ and $W$ is a weighting matrix. Thus, (5.1) corresponds to WLS with $W = V^{-1} = \operatorname{diag}(V_i^{-1}, 1 \leq i \leq n)$. This estimator is also known as the best linear unbiased estimator, or BLUE, denoted by $\hat{\beta}_{\text{BLUE}}$.

The IEE developed in previous sections can now be expressed more explicitly. Namely, if the $V_i$'s are known, one can calculate the BLUE by (5.1). On the other hand, if $\beta$ is known, the $V_i$'s can be estimated by MoM as follows [see (3.2)]:

$$(5.2) \qquad \hat{v}(j,k,l) = \frac{1}{n(j,k,l)} \sum_{i \in I(j,k,l)} (Y_{ij} - X_{ij}'\beta)(Y_{ik} - X_{ik}'\beta),$$

where $X_{ij}'$ is the $j$th row of $X_i$. When both $\beta$ and the $V_i$'s are unknown, one iterates between (5.1) and (5.2), starting with the OLS estimator. The latter is (5.1) with $V_i$ replaced by the identity matrix with the same dimension as $V_i$. The procedure is also called iterative reweighted least squares, or IRLS.

It is clear that IRLS is very easy to operate and is, in fact, no stranger to practitioners. However, even in this special case, very little is known about the convergence property. A procedure with some similarities has been used in robust linear regression (Huber [8]), and its convergence property has been studied. For example, Wolke and Schwetlick [18] considered a similar iterative procedure in solving robust regression problems involving an additional scale parameter. However, longitudinal data is certainly more complicated.

The results of Section 4 now have their versions in this special case.

COROLLARY 1. *Suppose that $\|V_i^{-1}\|$, $\|X_i\|$ and $\mathrm{E}(|Y_i|^4)$, $1 \leq i \leq n$, are bounded, Assumption A5 holds and $\liminf_{n \to \infty} \lambda_{\min}(n^{-1} \sum_{i=1}^n X_i' X_i) > 0$. Then we have $\mathrm{P}(\text{IRLS converges}) \to 1$, $\mathrm{P}[\sup_{m \geq 1} \{|\hat{\beta}^{(m)} - \hat{\beta}^*|/(p\eta)^{m/2}\} < \infty] \to 1$, and $\mathrm{P}[\sup_{m \geq 1} \{|\hat{v}^{(m)} - \hat{v}^*|/(R\eta)^{m/2}\} < \infty] \to 1$ as $n \to \infty$ for any $0 < \eta < (p \vee R)^{-1}$, where $(\hat{\beta}^*, \hat{v}^*)$ is the (limiting) IRLS estimator.*

COROLLARY 2. *Under the conditions of Corollary 1, $(\hat{\beta}^*, \hat{v}^*)$ is consistent.*

COROLLARY 3. *Under the same conditions of Corollary 1 except that Assumption A5 is strengthened to Assumption A5$'$, we have $\sqrt{n}(\hat{\beta}^* - \hat{\beta}_{\text{BLUE}}) \to 0$ in probability. Therefore, asymptotically, $\hat{\beta}^*$ is as efficient as the BLUE.*

NOTE. The condition $\liminf_{n \to \infty} \lambda_{\min}(n^{-1} \sum_{i=1}^n X_i' X_i) > 0$ in Corollary 1 needs some interpretation. Note that the term inside $(\cdots)$ is $n^{-1} X' X$,

where $X = (X_i)_{1 \leq i \leq n}$. Recall that the covariance matrix of the OLS estimator $\hat{\beta}_{\text{OLS}}$ is $(X'X)^{-1}(X'VX)(X'X)^{-1}$, where $V = \text{diag}(V_i, 1 \leq i \leq n)$. If the $V_i$'s are bounded from above and below, then we have $\text{Var}(\hat{\beta}_{\text{OLS}}) \sim (X'X)^{-1}$. Hence the condition is equivalent to the covariance matrix of $\sqrt{n}(\hat{\beta}_{\text{OLS}} - \beta)$ being bounded.

The proofs follow by verification of the assumptions of the theorems in Section 4. In the next section we study the empirical performance of IRLS.

## 6. Simulations.

6.1. *A simulation regarding Example 2.* In this subsection, we consider Example 2 in Section 2 (continued in Section 3) with simulated datasets. More specifically, we have $x'_{ij}\beta = \beta_0 + \beta_1 x_{ij}$, where the $x_{ij}$'s are generated from a $N(0, 1)$ distribution and then fixed throughout the simulation. We let $n = 100$. The $J_i$'s are specified as in Example 2, that is, $J_i = \{1, 3, 5\}$ or $\{2, 4\}$.

One potential advantage of IRLS is that it requires neither normality nor a parametric covariance model. In this simulation we consider three different scenarios. The first is when both normality and the parametric covariance model mentioned in Example 2 hold. The second is when normality fails but the parametric covariance model is still correct. In this scenario the true distribution of the random effects is *centralized exponential* [i.e., the distribution of $\sigma_u(\xi - 1)$, where $\xi \sim \text{Exponential}(1)$], while the distribution of $w$ and $e$ remains the same as in Scenario 1. The third is when both normality and the parametric covariance model no longer hold. In this scenario, the true distribution of $u$ and $w$ remains the same as in Scenario 2, while the true process of $w$ is MA(1) instead of AR(1) with normal noise $z$. In each scenario, we consider two cases of parameter values: (1) $\sigma_u^2 = 1$, $\sigma_w^2 = 9$, $\sigma_e^2 = 1$ and $\phi = 0.9$; (2) $\sigma_u^2 = 9$, $\sigma_w^2 = 25$, $\sigma_e^2 = 1$ and $\phi = 0.99$. Note that the second case corresponds to stronger within-subject correlations than the first case. In all cases, the true values for $\beta$ are $\beta_0 = 0.5$, $\beta_1 = 1.0$. These scenarios/cases will be denoted by 1.1, 1.2, and so forth.

We first study the convergence property of IRLS. The convergence criterion is $\max\{|\hat{\beta}_0^{(m)} - \hat{\beta}_0^{(m-1)}|, |\hat{\beta}_1^{(m)} - \hat{\beta}_1^{(m-1)}|\} + \max_{(j,k) \in D}\{|\hat{v}_{jk}^{(m)} - \hat{v}_{jk}^{(m-1)}|\} < 10^{-4}$. Note that here $d = |D| = 9$. In each case, we recorded the number of steps it took to converge. The results are summarized in Table 1.

Next, we compare the performance of IRLS with OLS and MLE, where MLE is that under normality and the parametric covariance model described in Example 2. Note that, under Scenarios 2 and 3, the MLE is not the true MLE. In Tables 2 and 3, we report the simulated means and covariance matrices as well as the true covariance matrix of the BLUE, even though the latter is not computable in practice. The covariance matrix of the BLUE serves as a lower bound for the covariance matrix of a WLS estimator. It should be pointed out that the IRLS estimator is not a WLS estimator and neither is the MLE, because both are nonlinear in $Y$. Nevertheless, we expect, by Corollary 3, the covariance matrix of the IRLS estimator to be close to that of the BLUE. The results for OLS and IRLS are based on 1,000 simulations. The results for MLE are a bit complicated. Although 1,000 simulations were run, not all resulted in convergence. Therefore, the results for MLE are based on the runs which converged. More specifically, for the Cases 1.1, 1.2, 2.1, 2.2, 3.1 and 3.2, the number of runs which converged for MLE were 798, 812, 794, 805, 749 and 724. (So, in the worst case, about 28% of the runs failed to converge.) In Table 2, $\hat{\beta}_{\text{OLS}}$, $\hat{\beta}_{\text{IRLS}}$ and $\hat{\beta}_{\text{MLE}}$ denote the simulated means of OLS, IRLS estimators and MLE. Similarly, in Table 3 $\hat{V}_{\text{OLS}}$, $\hat{V}_{\text{IRLS}}$ and $\hat{V}_{\text{MLE}}$ are the simulated covariance matrices and $V_{\text{BLUE}}$ is the true covariance matrix of the BLUE.

*Summary of results.* In most cases the IRLS algorithm converged in four to six steps. As for comparison of estimators, the main difference is in the

Table 1
*Number of steps to converge*

| Case | 2 | 3 | 4 | 5 | 6 | 7 | 8 | 9 | 10 | 11 |
|------|-----|-----|------|------|------|-----|-----|-----|-----|-----|
| 1.1  | 0.0 | 4.2 | 34.1 | 47.8 | 11.6 | 2.1 | 0.2 | 0.0 | 0.0 | 0.0 |
| 1.2  | 0.0 | 0.6 | 19.4 | 45.9 | 25.5 | 7.1 | 1.0 | 0.3 | 0.1 | 0.1 |
| 2.1  | 0.1 | 3.9 | 42.1 | 38.9 | 12.5 | 2.1 | 0.2 | 0.1 | 0.1 | 0.0 |
| 2.2  | 0.0 | 1.0 | 21.5 | 45.7 | 23.3 | 6.8 | 1.3 | 0.4 | 0.0 | 0.0 |
| 3.1  | 0.0 | 5.8 | 39.9 | 39.3 | 12.8 | 2.0 | 0.2 | 0.0 | 0.0 | 0.0 |
| 3.2  | 0.0 | 1.4 | 25.2 | 43.1 | 22.9 | 5.9 | 1.1 | 0.3 | 0.1 | 0.0 |

Numbers in the first row represent the steps; in each case, the numbers are percentages of times (out of a total of 1,000 simulations) that the IRLS converged after certain steps. Here Case 1.1 represents Scenario 1, Case 1, and so forth.

Table 2
*Simulated mean vectors*

| Case | Scenario 1 | | | Scenario 2 | | | Scenario 3 | | |
|------|-----------|-----------|-----------|-----------|-----------|-----------|-----------|-----------|-----------|
|      | $\hat{\beta}_{\text{OLS}}$ | $\hat{\beta}_{\text{IRLS}}$ | $\hat{\beta}_{\text{MLE}}$ | $\hat{\beta}_{\text{OLS}}$ | $\hat{\beta}_{\text{IRLS}}$ | $\hat{\beta}_{\text{MLE}}$ | $\hat{\beta}_{\text{OLS}}$ | $\hat{\beta}_{\text{IRLS}}$ | $\hat{\beta}_{\text{MLE}}$ |
| 1    | 0.50 | 0.50 | 0.50 | 0.50 | 0.49 | 0.50 | 0.51 | 0.51 | 0.52 |
|      | 0.99 | 1.00 | 1.00 | 1.00 | 1.00 | 1.00 | 1.01 | 1.01 | 1.01 |
| 2    | 0.50 | 0.51 | 0.50 | 0.51 | 0.51 | 0.53 | 0.50 | 0.47 | 0.56 |
|      | 1.00 | 1.00 | 1.00 | 1.00 | 1.01 | 1.00 | 1.02 | 1.02 | 1.01 |

In each case, the top row gives the simulated means for $\beta_0$; the bottom row gives those for $\beta_1$.

variances of the estimators of $\beta_1$. Both IRLS and MLE significantly outperform OLS. Furthermore, the performance of IRLS is very close to that of BLUE, although the latter is not computable in practice. Note that, when there is no correlation, OLS is the same as BLUE. As $\sigma_u^2$, $\sigma_w^2$ and $\phi$ increase, the within-subject correlation increases and, as a result, the difference between OLS and IRLS becomes larger. In the first two scenarios, MLE appeared to slightly outperform IRLS (and BLUE). This suggests that normality is not very important for the efficiency of the MLE for estimating $\beta$. However, the situation changes in Scenario 3 with IRLS (and BLUE) slightly outperforming MLE. In fact, the largest difference in favor of MLE is 18% less than IRLS in Case 2.2, while the largest difference in favor of IRLS is 14% less than MLE in Case 3.2. This suggests that the efficiency of MLE is likely to be undermined by misspecification of the covariance structure. It should be pointed out that, as mentioned above, the results for MLE are based on only runs which converged, the numbers of which are somewhere between 20 to 30 percent less. We expect the actual performance of MLE to be somewhat worse than reported here if all runs are reported, because the discarded runs may correspond to some "bad cases." Despite such concerns, the simulation results are consistent with our theoretical findings.

6.2. *A comparison with the one-step estimator.* In this subsection, we compare the small-sample performance of IEE with the one-step estimator $\hat{\beta}^{(3)}$ discussed in the second paragraph following Theorem 3. We consider the following simple example: $Y_{ij} = \beta_0 + \beta_1 x_i + e_{ij}$, $i = 1, \ldots, n$, $j = 1, 2$, where

TABLE 3
*Simulated covariance matrices*

| Case | $\hat{V}_{\text{OLS}}$ | | $\hat{V}_{\text{IRLS}}$ | | $\hat{V}_{\text{MLE}}$ | | $V_{\text{BLUE}}$ | |
|---|---|---|---|---|---|---|---|---|
| 1.1 | 9.80 | 0.11 | 10.36 | 0.36 | 10.00 | 0.48 | 9.30 | 0.11 |
|  | 0.11 | 5.51 | 0.36 | 2.09 | 0.48 | 1.87 | 0.11 | 2.16 |
| 1.2 | 36.41 | 0.61 | 36.41 | −0.10 | 35.57 | −0.08 | 34.10 | 0.06 |
|  | 0.61 | 17.41 | −0.10 | 1.30 | −0.08 | 1.23 | 0.06 | 1.29 |
| 2.1 | 9.87 | 0.14 | 9.94 | 0.07 | 9.19 | 0.02 | 9.30 | 0.11 |
|  | 0.14 | 5.35 | 0.07 | 2.25 | 0.02 | 2.08 | 0.11 | 2.16 |
| 2.2 | 38.74 | 1.78 | 38.48 | −0.01 | 34.54 | 0.15 | 34.10 | 0.06 |
|  | 1.78 | 15.86 | −0.01 | 1.31 | 0.15 | 1.11 | 0.06 | 1.29 |
| 3.1 | 7.14 | 0.44 | 7.25 | 0.41 | 6.77 | 0.45 | 7.15 | 0.23 |
|  | 0.44 | 5.34 | 0.41 | 4.01 | 0.45 | 4.13 | 0.23 | 3.74 |
| 3.2 | 25.48 | 0.26 | 26.22 | 0.22 | 24.73 | 0.46 | 25.40 | 0.63 |
|  | 0.26 | 17.55 | 0.22 | 10.80 | 0.46 | 12.35 | 0.63 | 9.61 |

In each case, the (2 × 2) simulated covariance matrices corresponding to OLS, IRLS and MLE and the true covariance matrix of the BLUE are presented. Here Case 1.1 represents Scenario 1, Case 1, and so forth.

the $x_i$'s are known covariates; $\beta_0$, $\beta_1$ are unknown regression coefficients; the $\varepsilon_{ij}$'s are independent with $\varepsilon_{ij} \sim N(0, \sigma_j^2)$, $j = 1, 2$.

Here we let $n = 10$. The $x_i$'s are generated from a Uniform$[0, 1]$ distribution, and then fixed throughout the simulation. The true parameters are chosen as $\beta_0 = 0.2$, $\beta_1 = 0.1$, $\sigma_1 = 1.0$ and $\sigma_2 = 4.0$. The results based on 1,000 simulations are reported in Table 4. It is seen that, once again, there is not much difference in terms of the simulated means. However, the simulated variances of the one-step estimator are about 31% and 28% larger than those of the IRLS estimator, which is the same as the MLE in this case, for the estimation of $\beta_0$ and $\beta_1$, respectively. Of course, both the one-step and IRLS estimators significantly outperformed the OLS estimator in terms of the simulated variance. It is remarkable that, although Theorem 1 (or Corollary 1) only ensures convergence with probability tending to one in large samples, with a fairly small sample size of $n = 10$, every one of our 1,000 IRLS runs actually converged. Here the criterion for convergence is similar to that in the previous subsection. The number of steps it took to converge ranged from 3 to 64 steps, and so are not reported here.

**7. A numerical example.** In this section we consider a data set presented by Hand and Crowder [5] regarding hip replacements of 30 patients. Each patient was measured four times, once before the operation and three times after, for hematocrit, TPP, vitamin E, vitamin A, urinary zinc, plasma zinc, hydroxyprolene (in milligrams), hydroxyprolene (index), ascorbic acid, carotine, calcium and plasma phosphate (twelve variables). One important feature of the data is that there is a considerable amount of missing observations. In fact, most of the patients have at least one missing observation for all twelve measured variables. In other words, the longitudinal data set is (seriously) unbalanced.

We consider the measured variable: hematocrit. This variable was considered by Hand and Crowder [5], who used the data to assess age, sex and time differences. The authors assumed an equicorrelated model and obtained Gaussian estimates of regression coefficients and variance components (i.e., MLE under normality). Here we take a robust approach without assuming a specific covariance structure. The covariates consist of the same variables as suggested by Hand and Crowder. They are an intercept, sex, occasion (three), sex by occasion interaction (three), age and age by sex interaction. For the hematocrit data the IRLS algorithm converged in seven steps. The results are shown in Table 5. The Gaussian (point) estimates of Hand and Crowder [5], page 106, are also included for comparison.

Table 4
*Comparison with the one-step estimator*

|  | Estimation of $\beta_0$ | | | Estimation of $\beta_1$ | | |
|---|---|---|---|---|---|---|
|  | OLS | One-step | IRLS | OLS | One-step | IRLS |
| Simulated mean | 0.15 | 0.17 | 0.18 | 0.16 | 0.16 | 0.15 |
| Simulated variance | 31.40 | 12.70 | 9.73 | 60.71 | 24.76 | 19.35 |

It is seen that the IRLS estimates are similar to the Gaussian estimates, especially for the parameters that are found significant. This is, of course, not surprising, because the Gaussian estimator and IRLS should both be close to the BLUE, provided that the covariance model suggested by Hand and Crowder is correct (the authors believed that their method was valid in this case). Taking into account the estimated standard errors, we found the coefficients $\beta_1$, $\beta_3$, $\beta_4$, $\beta_5$ and $\beta_6$ to be significant and the rest of the coefficients to be insignificant, where $\beta_1, \beta_2, \ldots$ are the coefficients corresponding to the variables described in the second paragraph of this section in that order. These are consistent with the findings of Hand and Crowder, with the only exception being $\beta_6$. Hand and Crowder considered jointly testing the hypothesis that $\beta_6 = \beta_7 = \beta_8 = 0$ and found an insignificant result. In our case, the coefficients are considered separately, and we found $\beta_7$ and $\beta_8$ to be insignificant and $\beta_6$ to be barely significant at 5% level. However, since Hand and Crowder did not publish the individual standard errors, this does not necessarily imply a difference.

**8. Discussion and concluding remarks.** In analysis of longitudinal data under linear mixed models, maximum likelihood (ML) and restricted maximum likelihood (REML) are well-established methods which apply to linear mixed models in general. What is the advantage of IRLS over these methods? First, the ML and REML methods require (i) normality; and (ii) correct specification of the parametric covariance structure. We do not believe that (i) is very important to the estimation of $\beta$, because the ML and REML estimators are consistent even without normality (Richardson and

TABLE 5
*Estimates for hematocrit*

| Coef. | $\beta_1$ | $\beta_2$ | $\beta_3$ | $\beta_4$ | $\beta_5$ | $\beta_6$ | $\beta_7$ | $\beta_8$ | $\beta_9$ | $\beta_{10}$ |
|---|---|---|---|---|---|---|---|---|---|---|
| IRLS | 3.19 | 0.08 | 0.65 | −0.34 | −0.21 | 0.12 | −0.051 | −0.051 | 0.033 | −0.001 |
| s.e. | 0.39 | 0.14 | 0.06 | 0.06 | 0.07 | 0.06 | 0.061 | 0.066 | 0.058 | 0.021 |
| Gaussian | 3.28 | 0.21 | 0.65 | −0.34 | −0.21 | 0.12 | −0.050 | −0.048 | 0.019 | −0.020 |

The first row gives IRLS estimates corresponding to, from left to right, intercept, sex, occasions (three), sex by occasion interaction (three), age and age by sex interaction; the second row gives estimated standard errors corresponding to the IRLS estimates; the third row gives the Gaussian maximum likelihood estimates that were obtained by Hand and Crowder [5].

Welsh [15], Jiang [9]). However, (ii) is crucially important. In fact, if (ii) fails, the ML and REML estimators may lose efficiency, even consistency (e.g., Liang and Zeger [12], White [17]). Note that, if the covariances are estimated using inconsistent estimators, the resulting "BLUE" is no longer BLUE, even asymptotically. In contrast, IRLS is not normality-based and does not require a parametric covariance structure. In this paper, we have shown that the IRLS estimator is asymptotically as efficient as the BLUE regardless of the covariance structure (and normality). This is confirmed by our simulation results. Second, IRLS is computationally easier to implement. The ML and REML methods require maximization of a nonlinear function or, at least, solution of a system of nonlinear equations. Although standard numerical procedures are available, problems and difficulties are often encountered in practice. For example, the Newton–Raphson procedure is known to be inefficient when the dimension of the solution is high, and its convergence is heavily affected by the choice of starting values; the EM algorithm is known to converge slowly. On the other hand, each step of IRLS is defined by closed-form expressions, therefore can be calculated analytically; the convergence is very fast, as we have demonstrated, and one does not need to worry about starting values. Note that because of the fast convergence of IRLS, in practice one may not need a lot of iterations (see the last paragraph of Section 4). This is also confirmed by our simulations. By the way, we have encountered difficulties in our simulations for computing the ML or quasi-likelihood estimator using the standard Newton–Raphson procedure. In fact, in some of our simulations, nearly thirty percent of the ML runs failed to converge. In contrast, every one of our IRLS runs converged and converged quickly, sometimes in two or three steps. Furthermore, the asymptotic properties of IEE, namely, Theorems 2 and 3, make sure that $\hat{v}^*$ always converges (in probability) to the true $v$, and hence results in an asymptotically optimal estimator of $\beta$. Other methods such as ML or REML do not necessarily ensure the latter property in case of misspecification of the variance–covariance structure. Therefore, IEE presents a robust and asymptotically efficient estimation procedure.

A basic assumption in this paper is Assumption A1, which essentially says that the number of different covariances between the observations is bounded. Although the assumption is required for establishing the theoretical properties of IEE, it does not mean that IEE cannot operate without such an assumption — it just makes it more difficult to justify. Alternatively, if the number of different covariances increases with the sample size, a parametric model for the covariance structure may be assumed to reduce the number of covariance parameters. Although, as mentioned earlier, such a model may be more sensitive to model misspecifications, it seems to be a reasonable approach in this difficult situation. For example, one may think of a linear model for the covariance components with factor covariates. Then, under Assumption A1, the number of columns of the design matrix in such a linear model is finite. More generally, one may incorporate continuous covariates in the linear model, which may be useful in some cases where the covariances of the observations depend on continuous covariates. Note that, however, the explicit expressions of the estimates of the covariance components that we derived in this paper may not exist under a parametric model. In fact, even under a linear model the estimates of the covariance components are subject to constraints, such as nonnegativity for the variance components, which are not guaranteed for the least-squares type estimates. Nevertheless, it is possible to develop a similar iterative procedure under a parametric covariance model.

IEE is not very picky about the initial estimator. For example, the starting point of IRLS is the OLS estimator; however, this is not an essential part of the algorithm. In fact, the only properties of OLS that are used are that (i) it is consistent, and (ii) $W = I$ is bounded from above and below (i.e., $\|W\|$ and $\|W^{-1}\|$ are bounded). Therefore, the starting point of IRLS may be replaced by another WLS estimator such that $W$ is bounded from above and below. A further question is whether one actually has *global* convergence, meaning that the algorithm converges regardless of the initial estimator (e.g., Luenberger [13]). It is an interesting question for which we do not have an answer at this point.

In this paper, we not only have proved that the IEE procedure converges under very mild conditions, we have shown that it converges linearly, that is, at an exponential rate. This is the main theoretical finding of this paper. Although the convergence of IEE would have been expected, the issue has never been rigorously addressed, especially regarding its convergence rate. Furthermore, the limiting IEE estimator is consistent and asymptotically efficient. These theoretical results are confirmed by our simulation studies. In addition, our method leads to consistent estimators of the covariances without assuming a parametric covariance structure and/or normality. Finally, we extended the robust estimation procedure (see Section 1) to unbalanced cases, which should make it a more attractive method for the analysis of longitudinal data.

**9. Proof of Theorem 1.** *Notation.* Throughout the rest, $\beta$, $v$, and so forth represent the true $\beta$, $v$, and so forth when there is no confusion. Given $\hat{v}$, the corresponding $\hat{\beta}$ is understood as $\hat{\beta} = \beta(\hat{v})$; similarly, given $\hat{\beta}$, the corresponding $\hat{v}$ is $\hat{v} = v(\hat{\beta})$ (see Section 4.1). Note that the two equations $\hat{\beta} = \beta(\hat{v})$ and $\hat{v} = v(\hat{\beta})$ need not hold simultaneously (unless at convergence). Therefore, for example, $v \circ \beta(\hat{v})$ may not equal $\hat{v}$. Also, we use notation such as $\beta^{(l)}$ for $\beta\{v^{(l)}\}$, and so forth, when there is no confusion.

The (spectral) norm of a matrix $A$ is defined as $\|A\| = \{\lambda_{\max}(A'A)\}^{1/2}$, and the 2-norm is defined as $\|A\|_2 = \{\text{tr}(A'A)\}^{1/2}$.

Let $\delta_0$ and $M_0$ denote positive numbers such that $\delta_0 \leq 1$, $M_0 \geq 1$ and $\delta_0 \leq M_0/16$. We define the following sets [see notation below (2.2), above (3.2) and above Theorem 1]:

$$B_0 = \{|\hat{\beta}^{(1)} - \beta| \leq \varepsilon_0\},$$

where $\varepsilon_0$ is a positive constant suitably chosen later on;

$$B_1 = \{\lambda_{\min}(\tilde{V}_i) \geq 2\delta_0 \text{ and } \lambda_{\max}(\tilde{V}_i) \leq M_0/2, 1 \leq i \leq n\},$$

where $\tilde{V}_i = (\tilde{v}_{ijk})_{j,k \in J_i}$ and $\tilde{v}_{ijk}$ is defined by the right-hand side of (3.2), if $v_{ijk} = v(j,k,l)$;

$$B_2 = \left\{ n(j,k,l)^{-1} \sum_{i \in I(j,k,l)} (1 + s_{ij})|Y_{ik} - \mu_{ik}| \leq \varepsilon_0^{-1/2}, 1 \leq l \leq L_{jk}, (j,k) \in D \right\};$$

$$B_3 = \left\{ \sup_{v \in \mathcal{V}} \|\partial \beta / \partial v\| \leq \varepsilon_0^{-1/4} \right\};$$

$B_4 = $ a set defined in Lemma 5 below such that $P(B_4) \to 1 \quad$ as $n \to \infty$;

$$B_5 = \left\{ n(j,k,l)^{-1} \left| \sum_{i \in I(j,k,l)} (\partial \mu_{ij}/\partial \beta_q)(Y_{ik} - \mu_{ik}) \right| \leq \varepsilon_0^{1/2}, \right.$$

$$\left. 1 \leq q \leq p, 1 \leq l \leq L_{jk}, (j,k) \in D \right\};$$

$C_0 = \{(1/4)\lambda_{\min}(V_i) \geq \delta_0 \text{ and } 4\lambda_{\max}(V_i) \leq M_0, 1 \leq i \leq n\};$

$C_1 = \{b^2(2\varepsilon_0^{1/2} + L_1^2\varepsilon_0^2) \leq \delta_0\};$

$C_2 = \{\sqrt{R(R \vee p)}[L_1^2 \varepsilon_0^{3/4} + (L_2 + 1)\varepsilon_0^{1/4}] \leq \eta/2\}.$

A rule for the notation is that $B$ represents a set involving both the $Y_i$'s and the $X_i$'s, while $C$ represents a set that only involves the $X_i$'s. Note that these are subsets of the probability space on which all the random variables are defined.

For the most part, the proof consists of two major steps. The first is to show that the sequences $\hat{\beta}^{(m)}$ and $\hat{v}^{(m)}$ ($m = 1, 2, \ldots$) are bounded. The second is to show that, within such a bounded range, the partial derivatives of $\tau(\cdot)$ and $\rho(\cdot)$ are bounded by some small numbers. We first state and prove some lemmas. Lemma 4 below ensures that, with high probability, $\hat{\beta}$ is confined to a neighborhood of $\beta$ implies that $\hat{v}$ is bounded; while Lemma 5 ensures the opposite, that is, with high probability, $\hat{v}$ is bounded implies that $\hat{\beta}$ is confined to the same neighborhood of $\beta$. Then, with Lemma 4 and Lemma 5, one is running a cycle, which allows one to argue that the entire sequence of $\hat{\beta}^{(m)}$, $\hat{v}^{(m)}$, $m = 1, 2, \ldots$, is bounded, starting from the beginning. Given that the sequence is bounded, the next thing is to obtain upper bounds for the derivatives, which is what Lemma 6 does. The proofs of Lemma 4 and Lemma 5 are fairly straightforward and therefore are omitted. Recall that $\varepsilon_1$ is the arbitrary positive constant involved in the definition of $L_j$, $j = 1, 2, 3$ (see Section 4).

LEMMA 4. *For any $0 < \varepsilon_0 \leq \varepsilon_1$, we have $\{|\hat{\beta} - \beta| \leq \varepsilon_0 \text{ implies } \lambda_{\min}(\hat{V}_i) \geq \delta_0 \text{ and } \lambda_{\max}(\hat{V}_i) \leq (9/16) M_0, \ 1 \leq i \leq n\} \supset B_1 \cap B_2 \cap C_1$.*

LEMMA 5. *Under Assumptions A3 and A4, for any $\varepsilon_0 > 0$, there is a set $B_4$ with $P(B_4) \to 1$ as $n \to \infty$ such that $\{\hat{v} \in \mathcal{V} \text{ implies } |\hat{\beta} - \beta| \leq \varepsilon_0\} \supset B_3 \cap B_4$.*

LEMMA 6. *For any $0 < \varepsilon_0 \leq \varepsilon_1$, we have*

$$\left\{ |\hat{\beta} - \beta| \leq \varepsilon_0 \text{ and } \hat{v} \in \mathcal{V} \text{ implies } \left| \frac{\partial \tau}{\partial \beta_q} \right|_{\beta = \hat{\beta}} \leq \eta, 1 \leq q \leq p, \text{ and} \right.$$

$$\left. \left| \frac{\partial \rho}{\partial v_r} \right|_{v = \hat{v}} \leq \eta, 1 \leq r \leq R \right\}$$

$$\supset B_2 \cap B_3 \cap B_5 \cap C_2.$$

PROOF. We have

(9.1) $$\frac{\partial \tau}{\partial \beta_q} = \sum_{r=1}^{R} \frac{\partial \beta}{\partial v_r} \cdot \frac{\partial v_r}{\partial \beta_q},$$

(9.2) $$\frac{\partial \rho}{\partial v_r} = \sum_{q=1}^{p} \frac{\partial v}{\partial \beta_q} \cdot \frac{\partial \beta_q}{\partial v_r}.$$

Let $v_r = \tilde{v}(j, k, l)$ for some $j$, $k$, $l$. Then, by (3.2), we have

$$\frac{\partial v_r}{\partial \beta_q} \bigg|_{\beta = \hat{\beta}} = -\frac{1}{n(j,k,l)} \sum_{i \in I(j,k,l)} \left\{ \frac{\partial g_j}{\partial \beta_q}(X_i, \hat{\beta}) \right\} \{Y_{ik} - g_k(X_i, \hat{\beta})\}$$

$$- \frac{1}{n(j,k,l)} \sum_{i \in I(j,k,l)} \left\{ \frac{\partial g_k}{\partial \beta_q}(X_i, \hat{\beta}) \right\} \{Y_{ij} - g_j(X_i, \hat{\beta})\}$$

$$= -(I_1 + I_2).$$

Furthermore, we have

$$I_1 = \frac{1}{n(j,k,l)} \sum_{i \in I(j,k,l)} \left\{ \frac{\partial g_j}{\partial \beta_q}(X_i, \beta) \right\} (Y_{ik} - \mu_{ik})$$

$$+ \frac{1}{n(j,k,l)} \sum_{i \in I(j,k,l)} \left\{ \frac{\partial g_j}{\partial \beta_q}(X_i, \hat{\beta}) - \frac{\partial g_j}{\partial \beta_q}(X_i, \beta) \right\}(Y_{ik} - \mu_{ik})$$

$$+ \frac{1}{n(j,k,l)} \sum_{i \in I(j,k,l)} \left\{ \frac{\partial g_j}{\partial \beta_q}(X_i, \hat{\beta}) \right\}\{g_k(X_i, \beta) - g_k(X_i, \hat{\beta})\}$$

$$= I_{11} + I_{12} + I_{13}.$$

Suppose that $|\hat{\beta} - \beta| \leq \varepsilon_0$. Then, it is easy to show that $|I_{13}| \leq L_1^2 \varepsilon_0$; $|I_{12}| \leq L_2 \varepsilon_0^{1/2}$ on $B_2$; and $|I_{11}| \leq \varepsilon_0^{1/2}$ on $B_5$. Thus, on $B_2 \cap B_5$, $|\hat{\beta} - \beta| \leq \varepsilon_0$ implies $|I_1| \leq L_1^2 \varepsilon_0 + (L_2 + 1)\varepsilon_0^{1/2}$, and, similarly, $|I_2| \leq L_1^2 \varepsilon_0 + (L_2 + 1)\varepsilon_0^{1/2}$. Therefore, by (9.1), we have, on $B_2 \cap B_3 \cap B_5 \cap C_2$, that $|\hat{\beta} - \beta| \leq \varepsilon_0$ and $v \in \mathcal{V}$ imply $|(\partial \tau / \partial \beta_q)|_{\beta = \hat{\beta}}| \leq \eta$. The arguments also show that, on $B_2 \cap B_5$, $|\hat{\beta} - \beta| \leq \varepsilon_0$ implies $|(\partial v / \partial \beta_q)|_{\beta = \hat{\beta}}| \leq 2\sqrt{R}\{L_1^2 \varepsilon_0 + (L_2 + 1)\varepsilon_0^{1/2}\}$. Thus, by (9.2), we have, on $B_2 \cap B_3 \cap B_5 \cap C_2$, that $|\hat{\beta} - \beta| \leq \varepsilon_0$ and $v \in \mathcal{V}$ imply $|(\partial \rho / \partial v_r)|_{v = \hat{v}}| \leq \eta$. □

We are now ready for the proof. For any $\delta > 0$, by Assumption A2, there is $M > 0$ such that $P(\mathcal{A}) > 1 - \delta$, where $\mathcal{A} = \{L_j \leq M, j = 1, 2, \text{ and } \|V_i\| \vee \|V_i^{-1}\| \leq M, 1 \leq i \leq n\}$. Then, choose $\delta_0$ and $M_0$ as in the third paragraph of this section such that $\lambda_{\min}(V_i) \geq 4\delta_0$ and $\lambda_{\max}(V_i) \leq M_0/4$, $1 \leq i \leq n$, on $\mathcal{A}$. Furthermore, if $\varepsilon_0 > 0$ is chosen such that

$$\varepsilon_0 \leq \varepsilon_1, \qquad b^2(2\varepsilon_0^{1/2} + M^2\varepsilon_0^2) \leq \delta_0,$$

$$\sqrt{R(R \vee p)}\{M^2\varepsilon_0^{3/4} + (M+1)\varepsilon_0^{1/4}\} \leq \frac{\eta}{2}$$

($\varepsilon_1$ is defined below Assumption A5), we have $C_j \supset \mathcal{A}$, $j = 0, 1, 2$.

By Assumption A3, there is $N_0$ such that $P(B_0^c) < \delta$, if $n \geq N_0$.

We have $\tilde{v}(j,k,l) - v(j,k,l) = n^{-1}(j,k,l) \sum_{i \in I(j,k,l)} \Delta_{i,j,k}$ for any $(j,k) \in D$, $1 \leq l \leq L_{jk}$, where $\Delta_{i,j,k} = (Y_{ij} - \mu_{ij})(Y_{ik} - \mu_{ik}) - v_{ijk}$. It follows, by Assumption A2, that $E\{\tilde{v}(j,k,l) - v(j,k,l)\}^2 = n^{-2}(j,k,l) \sum_{i \in I(j,k,l)} E(\Delta_{i,j,k}^2) \leq cn^{-1}(j,k,l)$ for some constant $c$. It is then easy to show, by Assumption A5, that there is $N_1$ such that $P(\max_i \|\tilde{V}_i - V_i\| > \delta_0) < \delta$, if $n \geq N_1$. Since $B_1 \supset \mathcal{A} \cap \{\max_i \|\tilde{V}_i - V_i\| \leq \delta_0\}$, we have $P(B_1^c) < 2\delta$, if $n \geq N_1$.

Similarly, there is $N_4$ such that $P(B_4^c) < \delta$, if $n \geq N_4$.

Now consider $B_2$. It is easy to see that $B_2 \supset \mathcal{A} \cap G$, where

$$G = \left\{ \frac{1}{n(j,k,l)} \sum_{i \in I(j,k,l)} |Y_{ij} - \mu_{ij}| \leq (M+1)^{-1}\varepsilon_0^{-1/2}, \right.$$

$$\left. 1 \leq l \leq L_{jk}, (j,k) \in D \right\}.$$

By Chebyshev's inequality and Assumptions A1, A2, we have

$$P(G^c) \le \sum_{1 \le l \le L_{jk}, (j,k) \in D} \frac{(M+1)\varepsilon_0^{1/2}}{n(j,k,l)} \sum_{i \in I(j,k,l)} E(|Y_{ij} - \mu_{ij}|)$$

$$\le c_1(M+1)\varepsilon_0^{1/2}$$

for some constant $c_1$. Thus, we have $P(B_2^c) < \delta + c_1(M+1)\varepsilon_0^{1/2}$.

Also, by Assumption A4, there is $\varepsilon_2 > 0$ such that $P(B_3^c) \le \delta$, if $\varepsilon_0 \le \varepsilon_2$.

Finally, for any $1 \le q \le p$, $(j,k) \in D$ and $1 \le l \le L_{jk}$, write

$$S_{q,j,k,l} = \left\{ \frac{1}{n(j,k,l)} \left| \sum_{i \in I(j,k,l)} \left( \frac{\partial \mu_{ij}}{\partial \beta_q} \right)(Y_{ik} - \mu_{ik}) \right| > \varepsilon_0^{1/2} \right\},$$

and $A(i) = \{|\partial \mu_{ij}/\partial \beta_q| \le M, v_{ikk} \le M\}$. Then we have $1_{A(i)} = 1 \ \forall i$ on $\mathcal{A}$. Thus, we have

$$P(S_{q,j,k,l} \cap \mathcal{A}) \le P\left\{ \frac{1}{n(j,k,l)} \left| \sum_{i \in I(j,k,l)} \left( \frac{\partial \mu_{ij}}{\partial \beta_q} \right)(Y_{ik} - \mu_{ik}) 1_{A(i)} \right| > \varepsilon_0^{1/2} \right\}$$

$$\le \frac{1}{\varepsilon_0 n^2(j,k,l)} \sum_{i \in I(j,k,l)} E\left\{ \left( \frac{\partial \mu_{ij}}{\partial \beta_q} \right)^2 v_{ikk} 1_{A(i)} \right\}$$

$$\le \frac{M^3}{\varepsilon_0 n(j,k,l)}.$$

Therefore, there is $N_5$ which depends on $M$ and $\varepsilon_0$, such that $P(S_{q,j,k,l} \cap \mathcal{A}) < \delta/h$, for all such $q$, $j$, $k$ and $l$, if $n \ge N_5$, where $h$ is the cardinality of the set $H = \{(q,j,k,l) : 1 \le q \le p, (j,k) \in D, l \le L_{jk}\}$, which is bounded by Assumption A1. It follows that

$$P(B_5^c \cap \mathcal{A}) \le \sum_{(q,j,k,l) \in H} P(S_{q,j,k,l} \cap \mathcal{A}) < \delta,$$

and hence $P(B_5^c) \le P(B_5^c \cap \mathcal{A}) + P(\mathcal{A}^c) < 2\delta$, if $n \ge N_5$.

In conclusion, for any $\delta > 0$, choose $\varepsilon_0 > 0$ such that $\varepsilon_0 \le \varepsilon_1 \wedge \varepsilon_2 \wedge 1$, $b^2(2\varepsilon_0^{1/2} + M^2\varepsilon_0^2) \le \delta_0$, $\sqrt{R(R \vee p)}\{M^2\varepsilon_0^{3/4} + (M+1)\varepsilon_0^{1/4}\} \le \eta/2$ and $c_1(M+1)\varepsilon_0^{1/2} < \delta$. Then, choose $N_j$, $j = 0, 1, 4, 5$, as above, and let $N$ be the maximum among them. It follows that $P(C_j^c) < \delta$, $j = 0, 1, 2$, and, when $n \ge N$, $P(B_j^c) < \delta$, $j = 0, 3, 4$, and $P(B_j^c) < 2\delta$, $j = 1, 2, 5$. Thus, with $\mathcal{E} = B_0 \cap \cdots \cap B_5 \cap C_0 \cap C_1 \cap C_2$, we have $P(\mathcal{E}^c) < 12\delta$, if $n \ge N$.

Now consider what happens on $\mathcal{E}$. For any $\beta^{(j)}$ such that $|\beta^{(j)} - \beta| \le \varepsilon_0$, $j = 1, 2$, we have

$$\tau\{\beta^{(1)}\} - \tau\{\beta^{(2)}\} = [\tau_q\{\beta^{(1)}\} - \tau_q\{\beta^{(2)}\}]_{1 \le q \le p}.$$

By Taylor expansion we have

$$\tau_q\{\beta^{(1)}\} - \tau_q\{\beta^{(2)}\} = \sum_{j=1}^{p} \frac{\partial \tau_q}{\partial \beta_j}\Big|_{\tilde{\beta}} \{\beta^{(1)} - \beta^{(2)}\},$$

where $\tilde{\beta}$ lies between $\beta^{(1)}$ and $\beta^{(2)}$ (but may depend on $q$). Note that $|\beta^{(j)} - \beta| \leq \varepsilon_0$, $j = 1, 2$, implies $|\tilde{\beta} - \beta| \leq \varepsilon_0$. Thus, by Lemma 4, we have $\tilde{v} \in \mathcal{V}$. Therefore, by Lemma 6 and the Cauchy–Schwarz inequality, we have

$$|\tau_q\{\beta^{(1)}\} - \tau_q\{\beta^{(2)}\}| \leq \sqrt{p}\eta|\beta^{(1)} - \beta^{(2)}|,$$

hence

(9.3) $$|\tau\{\beta^{(1)}\} - \tau\{\beta^{(2)}\}| \leq p\eta|\beta^{(1)} - \beta^{(2)}|.$$

Similarly, by Lemma 5 and Lemma 6 it can be shown that for any $v^{(j)} \in \mathcal{V}$, $j = 1, 2$, we have

(9.4) $$|\rho\{v^{(1)}\} - \rho\{v^{(2)}\}| \leq R\eta|v^{(1)} - v^{(2)}|.$$

Now, on $\mathcal{E}$ we have $|\hat{\beta}^{(1)} - \beta| \leq \varepsilon_0$ and $\hat{v}^{(1)} = \hat{v}^{(0)} \in \mathcal{V}$; $|\hat{\beta}^{(2)} - \beta| = |\hat{\beta}^{(1)} - \beta| \leq \varepsilon_0$, hence by Lemma 4, $\hat{v}^{(2)} \in \mathcal{V}$; thus, by Lemma 5, $|\hat{\beta}^{(3)} - \beta| \leq \varepsilon_0$, and $v^{(3)} = v^{(2)} \in \mathcal{V}$; and so on. It follows that $|\hat{\beta}^{(m)} - \beta| \leq \varepsilon_0$ and $\hat{v}^{(m)} \in \mathcal{V}$, $m = 1, 2, \ldots$.

We now apply (9.3) to obtain, for any $a, b > k$,

(9.5)
$$\begin{aligned}
|\hat{\beta}^{(2a-1)} - \hat{\beta}^{(2b-1)}| &= |\tau\{\hat{\beta}^{(2a-3)}\} - \tau\{\hat{\beta}^{(2b-3)}\}| \\
&\leq (p\eta)|\hat{\beta}^{(2a-3)} - \hat{\beta}^{(2b-3)}| \\
&\leq \cdots \\
&\leq (p\eta)^k|\hat{\beta}^{(2a-2k-1)} - \hat{\beta}^{(2b-2k-1)}| \\
&\leq 2\varepsilon_0(p\eta)^k.
\end{aligned}$$

It follows that the sequence $\{\hat{\beta}^{(2a-1)}, a = 1, 2, \ldots\}$ is a *Cauchy sequence* and hence convergent. Furthermore, since $\hat{\beta}^{(2a)} = \hat{\beta}^{(2a-1)}$, the entire sequence $\{\hat{\beta}^{(m)}, m = 1, 2, \ldots\}$ converges. By letting $b \to \infty$ and $k = a - 1$ in (9.5), it is easy to show that

$$\frac{|\hat{\beta}^{(m)} - \hat{\beta}^*|}{(p\eta)^{m/2}} \leq \frac{2\varepsilon_0}{p\eta} \equiv M_1, \qquad m = 1, 2, \ldots.$$

Similarly, by (9.4), it can be shown that $|\hat{v}^{(2a)} - \hat{v}^{(2b)}| \leq (R\eta)^k|\hat{v}^{(2a-2k)} - \hat{v}^{(2b-2k)}| \leq (2\sqrt{R}M_0)(R\eta)^k$ for any $a, b \geq k$. Thus, by similar arguments, the sequence $\{v^{(m)}, m = 1, 2, \ldots\}$ converges, and we have

$$\frac{|\hat{v}^{(m)} - \hat{v}^*|}{(R\eta)^{m/2}} \leq 2M_0\sqrt{R \vee \eta^{-1}} \equiv M_2, \qquad m = 1, 2, \ldots.$$

Let

$$\zeta_1 = \sup_{m \geq 1}\{|\hat{\beta}^{(m)} - \hat{\beta}^*|/(p\eta)^{m/2}\},$$

$$\zeta_2 = \sup_{m \geq 1}\{|\hat{v}^{(m)} - \hat{v}^*|/(R\eta)^{m/2}\}.$$

We have shown that for any $\delta > 0$, there is $N \geq 1$ such that, when $n \geq N$, $P(\text{IEE converges}) \geq P(\mathcal{E}) > 1 - 12\delta$ and $P(\zeta_j < \infty) \geq P(\zeta_j \leq M_j) \geq P(\mathcal{E}) > 1 - 12\delta$, $j = 1, 2$. This completes the proof.

J. JIANG  
DEPARTMENT OF STATISTICS  
UNIVERSITY OF CALIFORNIA  
ONE SHIELDS AVENUE  
DAVIS, CALIFORNIA 95616  
USA  
E-MAIL: jiang@wald.ucdavis.edu

Y. LUAN  
SCHOOL OF MATHEMATICS AND SYSTEM SCIENCES  
SHANDONG UNIVERSITY  
JINAN, SHANDONG 250100  
P. R. CHINA

Y.-G. WANG  
CSIRO MATHEMATICAL AND INFORMATION SCIENCES  
CSIRO LONG POCKET LABORATORIES  
120 MEIERS ROAD  
INDOOROOPILLY, QUEENSLAND 4068  
AUSTRALIA